\documentclass[11pt]{article} 
\usepackage{amsfonts,amsmath,latexsym,amssymb,mathrsfs,amsthm,comment}
\usepackage{graphicx}
\usepackage{xcolor}
\usepackage{booktabs}

\let\OLDthebibliography\thebibliography
\renewcommand\thebibliography[1]{
  \OLDthebibliography{#1}
  \setlength{\parskip}{1pt}
  \setlength{\itemsep}{0pt plus 0.0ex}
}

\def\numberlikeadb{\global\def\theequation{\thesection.\arabic{equation}}}
\numberlikeadb
\newtheorem{theorem}{Theorem}[section]
\newtheorem{lemma}[theorem]{Lemma}

\usepackage{lscape}
\usepackage{caption}
\usepackage{multirow}
\allowdisplaybreaks
\begin{document}

\title{On Stein factors in Stein's method for normal approximation
%On the solution of the standard normal Stein equation
}
\author{Robert E. Gaunt\footnote{Department of Mathematics, The University of Manchester, Oxford Road, Manchester M13 9PL, UK, robert.gaunt@manchester.ac.uk}}

\date{} 
\maketitle

\vspace{-5mm}

\begin{abstract} Building on the rather large literature concerning the regularity of the solution of the standard normal Stein equation, we provide a complete description of the best possible uniform bounds for the derivatives of the solution of the standard normal Stein equation when the test functions belong to the class of real-valued functions whose $k$-th order derivative is Lipschitz.
\end{abstract}

\noindent{{\bf{Keywords:}}} Stein's method; Stein equation; Stein factors; normal distribution

\noindent{{{\bf{AMS 2010 Subject Classification:}}} Primary 60F05}

%\section{Introduction}

\section{Introduction and main result}

Stein's method \cite{stein} is a powerful and widely-used technique for proving quantitative probabilistic limit theorems. Lying at the heart of Stein's method for normal approximation \cite{chen,np12,stein} is the  Stein equation
\begin{equation}\label{steineqn}
f'(x)-xf(x)=h(x)-Nh,    
\end{equation}
where $h:\mathbb{R}\rightarrow\mathbb{R}$ is a real-valued test function and $Nh$ denotes the expectation $\mathbb{E}[h(Z)]$ for $Z\sim N(0,1)$ a standard normal random variable. A solution to the Stein equation (\ref{steineqn}) is given by
\begin{equation}\label{steinsoln}
f_h(x)=\mathrm{e}^{x^2/2}\int_{-\infty}^x \{h(t)-Nh\}\mathrm{e}^{-t^2/2}\,\mathrm{d}t,  
\end{equation}
which is the unique bounded solution if $h$ is bounded or if $h$ is Lipschitz
%if $h^{(k)}$ is Lipschitz for $k\geq0$, where $h^{(0)}\equiv h$ 
(see, for example, \cite[Lemma 2.2]{ross}). In this note, we denote $\mathrm{Lip}_k(\mathbb{R})=\{h:\mathbb{R}\rightarrow\mathbb{R} \,|\,\text{$h^{(k)}$ is Lipschitz, $\|h^{(k+1)}\|_\infty$}<\infty\}$, $k\geq0$,  where $\|h^{(k+1)}\|_\infty$ denotes the
minimum Lipschitz constant of $h^{(k)}$, and $h^{(0)}\equiv h$.
%, where $\|h\|_{\mathrm{Lip}}=\sup_{x,y\in\mathbb{R}}|h(x)-h(y)|/|x-y|$. 
%Here $\mathrm{Lip}_0(\mathbb{R})=\mathrm{Lip}(\mathbb{R})$ is the class of Lipschitz functions on $\mathbb{R}$.

In deriving normal approximations via Stein's method, one typically requires bounds on the solution (\ref{steinsoln}) and its derivatives. Such bounds are often referred to as \emph{Stein factors}. Examples of such bounds are given by Lemma 2.4 of \cite{chen}: If $h:\mathbb{R}\rightarrow\mathbb{R}$ is bounded, then
\begin{equation}\label{first}
\|f_h\|\leq \sqrt{\pi/2}\|h(\cdot)-Nh\|, \quad \|f_h'\|\leq 2\|h(\cdot)-Nh\|,  
\end{equation}
whilst if $h:\mathbb{R}\rightarrow\mathbb{R}$ is Lipschitz, then
\begin{equation}\label{second}
\|f_h\|\leq 2\|h'\|, \quad \|f_h'\|\leq \sqrt{2/\pi}\|h'\|, \quad \|f_h''\|\leq 2\|h'\|,  
\end{equation}
where $\|\cdot\|=\|\cdot\|_\infty$ is the essential supremum norm on $\mathbb{R}$. For $g$ a Lipschitz function on $\mathbb{R}$, then $\|g'\|$ is its minimum Lipschitz constant. The bounds in (\ref{first}) and the third inequality of (\ref{second}) were proved earlier by \cite[Lemma II.3]{stein2}. The first bound in (\ref{second}) was later improved to $\|f_h\|\leq\|h'\|$ by \cite{d15}. The bounds in (\ref{first}) and (\ref{second}) are referred to as uniform bounds. Non-uniform bounds on $|f_h(x)|$ and $|f_h^{(k)}(x)|$, $k\geq1$, which are expressed in terms of the norms $\|h(\cdot)-Nh\|$ and $\|h^{(k)}\|$ but also have a dependence on the argument $x$, are given in, for example, \cite{chen,d15,es22,gaunt rate}. Henceforth, we will restrict our attention to uniform bounds.

In applications of Stein's method under which convergence rates faster than the classical $O(n^{-1/2})$ Berry-Esseen rate are sought (see, for example, \cite{bh85,daly,f18,gaunt rate,gr97,lu03}), the following generalisations of the inequalities of (\ref{second}) are often employed. Let $k\geq1$. Then, for $h\in \mathrm{Lip}_{k-1}(\mathbb{R})$,
%if $h^{(k-1)}:\mathbb{R}\rightarrow\mathbb{R}$ is Lipschitz, we have the bounds
\begin{equation}\label{third}
\|f_h^{(k-1)}\|\leq \frac{1}{k}\|h^{(k)}\|, \quad \|f_h^{(k)}\|\leq \frac{\Gamma((k+1)/2)}{\sqrt{2}\Gamma(k/2+1)}\|h^{(k)}\|, \quad \|f_h^{(k+1)}\|\leq 2\|h^{(k)}\|.   
\end{equation}
Note that the second and third inequalities of (\ref{third}) indeed reduce to the corresponding inequalities of (\ref{second}) in the case $k=1$, whilst the first inequality of (\ref{third}) reduces to the bound $\|f_h\|\leq\|h'\|$ of \cite{d15} when $k=1$. The first and second inequalities in (\ref{third}) were obtained by \cite{gr96} and \cite{gaunt rate}, respectively, using an alternative representation of \cite{b90,g91} for the solution $f_h$. The bound of \cite{gr96} was derived under an additional assumption that $h$ has three bounded derivatives. In the case $k=1$, this additional assumption was shown to not be required by \cite{d15}, whilst for $k\geq2$ this additional assumption can be seen to be unnecessary due to Proposition 2.1 of \cite{gms18}. The third bound of (\ref{third}) was obtained by \cite{daly}. Earlier, \cite[Lemma 5]{b86} had shown that $\|f_h^{(k+1)}\|\leq c_k\|h^{(k)}\|$ for some universal constant $c_k$ depending only on $k$.

%It is natural to ask whether further  

Motivated by the bounds in (\ref{third}), in this note we prove the following theorem.

\begin{comment}
which provides a complete description of the best possible bounds 
%of the form $\|f_h^{(j)}\|\leq c_k\|h^{(k)}\|$ 
uniform bounds 
for the solution (\ref{steinsoln}) of the standard normal Stein equation with test functions $h$ in the class $\mathrm{Lip}_k(\mathbb{R})$, $k\geq0$; that is the best possible bounds of the form $\|f_h^{(j)}\|\leq c_{j,k}\|h^{(k)}\|$ for $j,k\geq0$ and $c_{j,k}$ a universal constant depending only on $j$ and $k$.
\end{comment}

\begin{theorem}\label{thm1} Let $f_h$ denote the solution (\ref{steinsoln}) of the standard normal Stein equation (\ref{steineqn}), where $h\in\mathrm{Lip}_{k-1}(\mathbb{R})$, $k\geq1$. Then 

\begin{itemize}

\item[(i)] The bounds in (\ref{third}) are sharp.

\item[(ii)] Suppose that $k\geq1$ and $j\in\mathbb{Z}_{\geq0}\setminus\{k-1,k,k+1\}$.  Then, there does not exist a positive constant $C_{j,k}$, depending only on $j$ and $k$, such that the bound $\|f_h^{(j)}\|\leq C_{j,k}\|h^{(k)}\|$ holds for all $h\in\mathrm{Lip}_{k-1}(\mathbb{R})$.
%, $k\geq0$.

\end{itemize}
\end{theorem}

%Theorem \ref{thm1} tells us that amongst

\begin{comment}
Theorem \ref{thm1} provides a complete description of the best possible bounds 
%of the form $\|f_h^{(j)}\|\leq c_k\|h^{(k)}\|$ 
uniform bounds 
for the solution (\ref{steinsoln}) of the standard normal Stein equation with test functions $h$ in the class $\mathrm{Lip}_k(\mathbb{R})$, $k\geq0$; that is the best possible bounds of the form $\|f_h^{(j)}\|\leq c_{j,k}\|h^{(k)}\|$ for $j,k\geq0$ and $c_{j,k}$ a universal constant depending only on $j$ and $k$. 
%In particular, we conclude that the bounds of (\ref{third}) are sharp and that
\end{comment}

It is already known that the first and third inequalities of (\ref{third}) are sharp (see \cite[p.\ 568]{daly}), although our verification that the second inequality is sharp is new. Part (ii) of Theorem \ref{thm1} is of interest because in applications of Stein's method for normal approximation it is generally preferable to impose weaker differentiability conditions on the test functions $h$; for example, this allows quantitative probabilistic limit theorems to be stated in stronger probability metrics. This result puts constraints on how much we can weaken the differentiability conditions on the test functions $h$.
%We thus conclude
%us that if we require bounds on the supremum norm $\|f_h^{(k)}\|$ then we are required to take $h\in\mathrm{Lip}_{k-1}(\mathbb{R})$
%us that from this perspective we cannot in general do better than the bound $\|f_h^{(k+1)}\|\leq 2\|h^{(k)}\|$ of \cite{daly} when working with test function $h$ from the class $\mathrm{Lip}_k(\mathbb{R})$. 
We remark that part (ii) of Theorem \ref{thm1} was proven for $j\geq k+3$ 
%$j\in\{0,1,\ldots,k-3\}$ 
by \cite{gl23}. Our contribution is to extend this result to $j\in\mathbb{Z}_{\geq0}\setminus\{k-1,k,k+1\}$ from which we conclude that when working with test functions $h$ from the class $\mathrm{Lip}_{k-1}(\mathbb{R})$ the only (sharp) uniform bounds of the form $\|f_h^{(j)}\|\leq c_{j,k}\|h^{(k)}\|$, for $j\geq0$, $k\geq1$ and $c_{j,k}$ a universal constant depending only on $j$ and $k$, that hold for all $h\in\mathrm{Lip}_{k-1}(\mathbb{R})$ are given by the inequalities in (\ref{third}).
%(which are also sharp).

%Indeed, Theorem \ref{}

%We provide a simple verification that the second inequality of (\ref{third}) is sharp in our proof of Theorem \ref{thm1}. The more sub

%\section{Main result and proof}

\section{Proof of Theorem \ref{thm1}}

We will need the following lemma, which is stated in \cite[Lemma 3.3]{gl23} and can also be read off as an intermediate bound in the proof of Theorem 3.1 of \cite{gaunt rate}.  

\begin{lemma}\label{lem1}[Gaunt \cite{gaunt rate}]
Let $X_1,\ldots,X_n$ be independent and identically distributed random variables such that $\mathbb{E}[|X_1|^{p+1}]<\infty,$ $p \geq 2.$ Suppose further that $\mathbb{E}[X_1^k] =\mathbb{E}[Z^k] $ for $1 \leq k \leq p$, where $Z\sim N(0,1)$. Write $W_n =n^{-1/2}\sum^n_{i=1} X_i.$ Let $f_h$ be the solution (\ref{steinsoln}) of the standard normal Stein equation (\ref{steineqn}).
%, for $h\in\mathrm{Lip}_{p-1}(\mathbb{R})$. 
Then 
\begin{equation}\label{iid example}
|\mathbb{E}[h(W_n)]-\mathbb{E}[h(Z)]| \leq \frac{\|f_h^{(p)} \|}{n^{(p-1)/2}}\bigg\{\frac{\mathbb{E}[|X_1|^{p-1}]}{(p-1)!}+\frac{\mathbb{E}[|X_1|^{p+1}]}{p!}\bigg\}.
\end{equation}
\end{lemma}

\noindent{\emph{Proof of Theorem \ref{thm1}.}} (i) As noted in the Introduction, the first and third inequalities of (\ref{third}) are already known to be sharp, so we consider only the second inequality. We begin by recalling equation (2.6) of \cite{gaunt rate} specialised to the standard normal setting: for $x\in\mathbb{R}$,
\begin{equation}\label{fourth}
f_h^{(k)}(x)=-\int_0^\infty \frac{\mathrm{e}^{-(k+1)s}}{\sqrt{1-\mathrm{e}^{-2s}}}\mathbb{E}\big[Zh^{(k)}(\mathrm{e}^{-s}x+\sqrt{1-\mathrm{e}^{-2s}}Z)\big]\,\mathrm{d}s,   
\end{equation}
where $Z\sim N(0,1)$. For $a>0$, we now suppose that $h_a^{(k)}$ is a continuous function such that $h_{a}^{(k)}(x)=1$ for $x\geq a$, $h_{a}^{(k)}(x)=-1$ for $x\leq -a$, and $h_a^{(k)}(x)$ is linear for $x\in(-a,a)$. Note that $h_a\in \mathrm{Lip}_{k-1}(\mathbb{R})$ with $\|h_a^{(k)}\|=1$. Then, as $a\rightarrow0$,
\begin{align*}
f_{h_a}^{(k)}(0)&=-\int_0^\infty \frac{\mathrm{e}^{-(k+1)s}}{\sqrt{1-\mathrm{e}^{-2s}}}\mathbb{E}\big[Zh_a^{(k)}(\sqrt{1-\mathrm{e}^{-2s}}Z)\big]\,\mathrm{d}s  \\
&\rightarrow-\int_0^\infty \frac{\mathrm{e}^{-(k+1)s}}{\sqrt{1-\mathrm{e}^{-2s}}}\mathbb{E}[|Z|]\,\mathrm{d}s\\
&=-\sqrt{\frac{2}{\pi}}\int_0^\infty \frac{\mathrm{e}^{-(k+1)s}}{\sqrt{1-\mathrm{e}^{-2s}}}\,\mathrm{d}s=-\frac{\Gamma((k+1)/2)}{\sqrt{2}\Gamma(k/2+1)},
\end{align*}
where we evaluated the final integral by using equation (2.7) of \cite{gaunt rate}. We have therefore shown that the second inequality of (\ref{third}) is sharp.

%On taking $h(x)=x^k/k!$ (which is in the class $\mathrm{Lip}_{k-1}(\mathbb{R})$ with $\|h^{(k)}\|=1$) in (\ref{fourth}) and evaluating the integral $\int_0^\infty \frac{\mathrm{e}^{-(k+1)s}}{\sqrt{1-\mathrm{e}^{-2s}}}\,\mathrm{d}s$ (see \cite[equation (2.7)]{gaunt rate}) we readily deduce that the second inequality of (\ref{third}) is sharp.

\vspace{2mm}

\noindent (ii) We first prove that a bound of the form $\|f_h^{(j)}\|\leq C_{j,k}\|h^{(k)}\|$ cannot hold for all $h\in\mathrm{Lip}_{k-1}(\mathbb{R})$ if $j=0,1,\ldots,k-2$, with $k\geq2$. Suppose to the contrary that there exists a bound of the form $\|f_h^{(p)}\|\leq C_{p,k}\|h^{(k)}\|$ that holds for all $h\in\mathrm{Lip}_{k-1}(\mathbb{R})$ if $p=0,1,\ldots,k-2$, with $k\geq2$. In that case, for random variables $X_1,\ldots,X_n$ and $W_n$ defined as in Lemma \ref{lem1}, we would obtain by an application of inequality (\ref{iid example}) that, for $h\in\mathrm{Lip}_{k-1}(\mathbb{R})$,
\begin{equation}\label{fifth}
|\mathbb{E}[h(W_n)]-\mathbb{E}[h(Z)]| \leq \frac{C_{p,k}\|h^{(k)} \|}{n^{(p-1)/2}}\bigg\{\frac{\mathbb{E}[|X_1|^{p-1}]}{(p-1)!}+\frac{\mathbb{E}[|X_1|^{p+1}]}{p!}\bigg\}.
\end{equation}
Consider now the test function $h_p(x)=x^{p+1}$, which is in the class $\mathrm{Lip}_{k-1}(\mathbb{R})$, since $p\leq k-2$. Moreover, $h_p^{(k)}(x)=0$, and so from (\ref{fifth}) we deduce that $\mathbb{E}[h_p(W_n)]-\mathbb{E}[h_p(Z)]=0$, that is $\mathbb{E}[W_n^{p+1}]=\mathbb{E}[Z^{p+1}]$. However, according to \cite[Lemma 3.5]{gaunt normal}, $\mathbb{E}[W_n^{p+1}]=\mathbb{E}[Z^{p+1}]+(\mathbb{E}[X_1^{p+1}]-\mathbb{E}[Z^{p+1}])/n^{(p-1)/2}$, and so taking $X_1$ to be such that $\mathbb{E}[X_1^{p+1}]\not=\mathbb{E}[Z^{p+1}]$ leads to a contradiction.

Now we prove that a bound of the form $\|f_h^{(j)}\|\leq C_{j,k}\|h^{(k)}\|$ cannot hold for all $h\in\mathrm{Lip}_{k-1}(\mathbb{R})$ if $j=k+2,k+3,\ldots$, with $k\geq1$. We shall prove this by selecting a function $h\in\mathrm{Lip}_{k-1}(\mathbb{R})$ for which $f_h^{(k+1)}(x)$ is not differentiable.
%It suffice to consider the case $j=k+2$, since \cite{gl23} have already proved the result for $j\geq k+3$. 
We begin by recalling the following representation of \cite[Lemma 2.2]{daly} for $f_h^{(k+1)}(x)$; for $x\in\mathbb{R}$,
\begin{equation*}
f_h^{(k+1)}(x)=h^{(k)}(x)-\frac{Z^{(k+1)}(-x)}{(k-1)!}A_k(x)-\frac{Z^{(k+1)}(x)}{(k-1)!}B_k(x),   
\end{equation*}
where $Z(x):=\Phi(x)/\phi(x)$ with $\phi$ and $\Phi$ denoting the standard normal density and cumulative distribution function, respectively, and
\begin{align*}
A_k(x):=\int_{-\infty}^x h^{(k)}(t)\phi(t)Z^{(k-1)}(t)\,\mathrm{d}t, \quad
B_k(x):=\int_x^\infty h^{(k)}(t)\phi(t)Z^{(k-1)}(-t)\,\mathrm{d}t.
\end{align*}
We now consider a test function $h:\mathbb{R}\rightarrow\mathbb{R}$ which is such that $h^{(k-1)}(x)=|x|$, and so $h\in\mathrm{Lip}_{k-1}(\mathbb{R})$. We then have $h^{(k)}(x)=1$ for $x>0$ and $h^{(k)}(x)=-1$ for $x<0$. The function $Z^{(k+1)}(x)$ is continuous at $x=0$ (this is easily seen from part i of Lemma 2.1 of \cite{daly}) and $A_k(x)$ and $B_k(x)$ are also continuous at $x=0$ due to the continuity of integrals of Lebesgue-integrable functions. Hence, $f_h^{(k+1)}(0+)-f_h^{(k+1)}(0-)=h^{(k)}(0+)-h^{(k)}(0-)=2$, and so $f_h^{(k+1)}(x)$ is not continuous at $x=0$, and thus not differentiable at $x=0$. This completes the proof. \hfill\qed

\section*{Acknowledgements}
The author is funded in part by EPSRC grant EP/Y008650/1 and EPSRC grant UKRI068.

\end{document}